\documentclass[12pt]{amsart}
\usepackage{amsmath,amsthm,latexsym,amscd,amsbsy,amssymb,url}
\usepackage[all]{xypic}
 \relax

\chardef\bslash=`\\ % p. 424, TeXbook

\makeatletter
\def\verbatim{\interlinepenalty\@M \@verbatim
  \leftskip\@totalleftmargin\advance\leftskip2pc
  \frenchspacing\@vobeyspaces \@xverbatim}
\makeatother
\newtheorem{thm}{Theorem}[section]
\newtheorem{cor}[thm]{Corollary}
\newtheorem{lem}[thm]{Lemma}
\newtheorem{pro}[thm]{Proposition}

\newtheorem{que}[thm]{Question}

\begin{document}

%%%%%%% Begin Topmatter %%%%%%%%%%

\title
{On Topological Groups of Automorphisms}
\author{Raushan  Buzyakova}
\email{Raushan\_Buzyakova@yahoo.com}


%%%%%%% End topmatter %%%%%%%%%

\begin{abstract}{
We study spaces $X$ for which the space $Hom_p(X)$ of automorphisms with the topology of point-wise convergence is a topological group. We identify large classes of spaces $X$  for which $Hom_p(X)$ is or is not  a topological group. 
}
\end{abstract}

\maketitle
\markboth{R. Buzyakova}{On Topological Groups of Automorphisms}
{ }

\section{Introduction}\label{S:introduction}
\par\bigskip\noindent
The focus of this study is the set of homeomorphic bijections (automorphisms) of a space $X$ with itself, denoted by $Hom(X)$, endowed with the topology of point-wise convergence. The resulting space is denoted by $Hom_p(X)$. Recall that a standard basic set of this topology is in the form $\{h\in Hom(X): h(x_i)\in O_i, i=1,...,n\}$, where $x_1,...,x_n$ are arbitrary fixed elements of $X$ and $O_1, ...,O_n$ are arbitrary fixed open subsets of $X$. It is known that $Hom(X)$ is an algebraic group with respect to the operation of composition, which is always assumed to be the group operation in this study.  It is well-known that $Hom_p(X)$ need not be a topological group. For example,  $Hom_p(\mathbb R^2)$ is not a topological group  as mentioned in one of the exercises in \cite{AT}. There are also many examples of spaces for which the respective structure is a topological group. For example, $Hom_p(L)$ is a topological group for any connected linearly ordered space $L$ (\cite{Sor}). Moreover, it was shown in \cite{B} that $Hom_p(X)$ is a topological group whenever $X$ is the union of finitely many closed connected  linearly ordered subspaces. In this study, will will focus on finding spaces $X$ for which $Hom_p(X)$ is a  topological group. To justify a clear lack of homogeneity in our candidates, we first start by proving that too much of homogeneity in $X$ is often the reason for $Hom_p(X)$ not being a topological group. Clearly,  $Hom_p(X)$ is a topological group if $X$ has no automorphisms except the identity. The spaces in this study are rich with automorphisms.

\par\smallskip\noindent
In notation and terminology we will follow \cite{Eng}. For general facts on topological groups, we refer the reader to \cite{AT}. To distinguish ordered pairs from open intervals, the former will be denoted by $\langle a, b\rangle$ and the latter by $(a,b)$. All spaces are assumed Tychonov. If $S$ and $T$ are subsets of $Hom(X)$ for some $X$, then $S\circ T$ and $S^{-1}$ are the sets $\{f\circ g: f\in S, g\in T\}$ and $\{f^{-1}: f\in S\}$, respectively.

\par\bigskip\noindent
\section{Ones That Are Not}

\par\bigskip\noindent
Let us start by identifying a large class of spaces that cannot be a source for our examples. Among those are spaces with high degree of homogeneity, namely, the spaces that are strongly $n$-homogeneous for each positive integer $n$ and have no isolated points. Recall that a space $X$ is {\it strongly $n$-homogeneous} if given any two $n$-sized subsets $A$ and $B$ of $X$ and any bijection $h:A\to B$ there exists $f\in Hom(X)$ such that $f|{_A} = h$. This concept is as classical and well-studied as the concept of homogeneity itself. We would like to reference one paper in connection with this concept, namely, that of G. Ungar \cite{U} due to the fact that the author of that paper studies  $Hom(X)$ as well but with the compact-open topology.
 Many classical objects are strongly $n$-homogeneous. In particular, the space of rationals, the space of irrationals, the Cantor Set, and $\mathbb R^n$ for any $n>1$ have the property. Note that the space of reals is not such since  any homeomorphism on the reals  must preserve or reverse the order.

\par\bigskip\noindent
The arguments of Lemmas \ref{lem:not1} and \ref{lem:not11} are very standard and have been used by many when exercising in proofs that $Hom_p(\mathbb R^2)$ or alike is not a topological group. The goal of discussing these folklore results in the presented forms is to justify the routes we chose in the next section in our search for spaces $X$ for which $Hom_p(X)$ is a topological group.  

\par\bigskip\noindent
\begin{lem}\label{lem:not1} Let $X$ have no isolated points. If $X$ is strongly $n$-homogeneous for each $n$, then the map $\langle f, g\rangle\mapsto f\circ g$ is not continuous on $Hom_p(X)\times Hom_p(X)$.
\end{lem}
\begin{proof}
Fix arbitrary $f,g\in Hom(X)$, $x\in X$, and an open neighborhood $O_z$ of $z=f(g(x))$ such that $X\setminus \overline O_z\not = \emptyset$. Let $V_{f\circ g} = \{h\in Hom(X): h(x)\in O_z\}$. Let $V_f$ and $V_g$ be arbitrary neighborhoods of $f$ and $g$, respectively. Our goal is to show that  $V_f\circ V_g$ is not a subset of $V_{f\circ g}$. 

Put  $y=g(x)$. We may assume that there exist $y_1=y, y_2,..., y_n$ and open $O_1, ..., O_n$ such that $V_f= \{h\in Hom(X): h(y_i)\in O_i, i=1,...,n\}$. 
Similarly, we may assume that there exist $x_1=x, x_2,..., x_m$ and open $B_1, ..., B_m$ such that $V_g= \{h\in Hom(X): h(x_i)\in B_i, i=1,...,m\}$. 
Since $X$ has no isolated points, we can pick $x^*\in X \setminus \{x, x_2, ...., x_m\}$ and $y^*\in B_1\setminus \{g(x_i): i = 1,...,m\}$.  By strong $m$-homogeneity, there exists $g_1\in Hom(X)$ such that $g_1(x^*)=y, g_1(x)=y^*$, and $g_1(x_i)=g(x_i)$ for $i=2,...,m$. Clearly $g_1\in V_g$. Next pick an arbitrary $z^*\in X\setminus O_z$.  By strong $n$-homogeneity, there exists $f_1\in Hom(X)$ such that $f_1(y^*)=z^*$, and $f_1(y_i)=f(x_i)$ for $i=1,...,n$. Clearly $f_1\in V_f$. Since $f_1(g_1(x))= z^*\not \in O_z$, our goal is achieved.
\end{proof}

\par\bigskip\noindent
\begin{lem}\label{lem:not11} Let $X$ have no isolated points. If $X$ is strongly $n$-homogeneous for each $n$, then the map $f\mapsto f^{-1}$ is not continuous on $Hom_p(X)$.
\end{lem}
\begin{proof}
Fix an arbitrary $f\in Hom(X)$, $x\in X$, and an open neighborhood $O_x$ of $x$ such that $X\setminus \overline O_x\not = \emptyset$. Let $V = \{h\in Hom(X): h(y)\in O_x\}$, where $f(x)=y$. Then $V $ is an open neighborhood  of $f^{-1}$. Let $U$ be an arbitrary neighborhood of $f$. It suffices to show that $U^{-1}$ is not a subset of $V$. We may assume that there exist a finite $n$-sized subset  $\{x_i: i=1,...,n\}$
of $X$ and a collection of open sets $\{O_i:i=1,...,n\}$ of $X$ such that $U=\{h\in Hom(X): h(x_i)\in O_i\ for\ i=1,...,n\}$. Next, fix an arbitrary $z$ in $X\setminus (O_x\cup \{x_i:i=0,...., n-1\})$. Such a $z$ exists since $X$ has no isolated points and $O_x$ is small enough. Since $X$ is strongly $(n+1)$-homogeneous and has no isolated points, there exists $g\in Hom(X)$ such that $g(z)=y$ and $g(x_i)\in O_i\setminus \{y\}$. Then $g\in U$. Since $g^{-1}(y)=z\not \in O_x$, we conclude that $g^{-1}\not\in V$.
\end{proof}

\par\bigskip\noindent
\begin{cor}\label{cor:not111}
If $X$ is strongly $n$-homogeneous for each $n$  and has no isolated points, then $Hom_p(X)$ is neither topological, nor paratopological, nor semitopological group.
\end{cor}

\par\bigskip\noindent
In Corollary \ref{cor:not111} and the preceding lemmas, we cannot replace "strong $n$-homogeneity" with "$n$-homogeneity".  Recall that $X$ is $n$-homogeneous if given any two $n$-sized subsets  $A$ and $B$ of $X$ there exists $h\in Hom(X)$ such that $f(A) = B$ (see \cite{U} for references). Note that $\mathbb R$ is $n$-homogeneous for any positive $n$ but not strongly $n$-homogeneous for $n>2$. Therefore, $\mathbb R$ is a witness that "{\it strongly}" cannot simply be dropped from the mentioned statements.  

\par\bigskip\noindent
For the remainder of this study, we say that $X\subset Y$ is {\it $h$-embedded }in $Y$ if any homeomorphisim $h$ from $X$ onto $X$ can be extended to a homeomorphism from $Y$ onto $Y$.

\par\bigskip\noindent
\begin{que}\label{que:not2}
Suppose that $X$ is $h$-embedded in $Y$ as an open subset.
Suppose that $Hom_p(X)$ is not a topological group. Can one conclude that  $Hom_p(Y)$ is not a topological group? What if $X$ is a dense (or open and dense) subset of $X$? 
\end{que}

\par\bigskip\noindent
In Question \ref{que:not2},  it is important to place a strong condition on how  $X$ is $h$-embedded in $Y$. This is justified by the following statement.

\par\bigskip\noindent
\begin{pro}\label{pro:not3}
For any  space $X$  there exists a space $Y$ such that $X$ $h$-embeds in $Y$ as a closed subspace and  $Hom_p(Y)$ is a topological group.
\end{pro}
\begin{proof}
The conclusion is a corollary to Corollary \ref{cor:yes000}. The proof is in the next section since it has a positive flavor and does not match the goal of this section.

\end{proof}

\par\bigskip
\section{Ones That Are}

\par\bigskip\noindent
We start by showing that for any space $X$ there exists a space $Y$ such that $X$  $h$-embeds in a $Y$ as a closed subspace and $Hom_p(Y)$ is a topological group. To construct such  $Y$ we introduce a structure that looks very similar to the Alexandroff double but instead of an isolated twin $x'$ for each $x$ we attach to each $x$ a segment. Next is a detailed construction.

\par\bigskip\noindent
{\bf Construction of Alexandroff Noodles.} Given a topological space $X$, the Alexandroff Noodles of $X$, denoted as $XI$, is defined as follows. 
Let $D=\{d_x:x\in X\}$ be a discrete space of the same cardinality as $X$ indexed by the  elements of $X$. Then $XI$ is the quotient space defined by the partition on $X\cup (D\times [0,1])$ whose only non-trivial elements are in the form $\{x,\langle d_x,0\rangle\}$, where $x\in X$. 
For brevity, we agree to refer to elements $\{x,\langle d_x,0\rangle\}$ by $0_x$, to $\{ \langle d_x, r\rangle \}$ for $r\in (0,1]$ by $r_x$ and to 
$\{d_x\}\times (a,b)$ by $(a_x, b_x)$. That is, every  point of $XI$ is identified as $r_x$ with $r\in [0,1]$ and $x\in X$.

\par\bigskip\noindent
Clearly, if $X$ is Tychonoff, so is $XI$.  
Note that $O_X= \{0_x:x\in X\}$  is a closed subspace of $XI$ and is homeomorphic to $X$ by virtue of map $x\mapsto 0_x$. Further,  let $f$ be a homeomorphism from $O_X$ onto $O_X$. The map $g: XI\to XI$ defined by $g(r_x) = r_y$, where $f(0_x)=0_y$ is a homeomorphic extension of $f$. Thus, $X$ is $h$-embedded in $XI$ as a closed subspace.
We will next show that $Hom_p(XI)$ is a topological group for any $X$.

\par\bigskip\noindent
\begin{lem}\label{lem:yes0}
Let $X$ be a topological space and let $XI$ be its Alexandroff Noodles. Then, the map $f\mapsto f^{-1}$ is continuous on $Hom_p(XI)$.
\end{lem}
\begin{proof}
Fix $f\in Hom(XI)$, $p\in XI$, and an open neighborhood $O_p$ of $p$. Put $q = f(p)$ and $V = \{h\in Hom(XI): h(q)\in O_p\}$. Clearly, $V$ is an open neighborhood of $f^{-1}$. It remains to find an open neighborhood $U$ of $f$ such  that $U^{-1}\subset V$. We will break down our argument into cases.
\begin{description}
	\item[\it Case \rm ($p=0_x$ for some $x\in X$)] If $x$ is not isolated in $X$, then $q = 0_y$ for some $y\in X$.
 Put $U= \{h\in Hom(XI): h(1_x)\in (0.5_y, 1_y]\}$. Clearly, $f\in U$. If $h\in U$, then $h(p) = 0_y$. Hence, $U^{-1}\subset V$.
\par\smallskip\noindent
If $x$ is isolated in $X$, then $q=0_y$ or $q=1_y$ for an isolated $y\in X$. In either case, $[0_y, 1_y]$ is open in $XI$. Assume $q=0_y$.
Then  set $U= \{h\in Hom(XI): h(p)\in [0_y, 0.5_y)\}$ is as desired. Indeed, if $h\in U$, then $h(p) = 0_y$. Hence, $U^{-1}\subset V$.
	\item[\it Case \rm ($p=1_x$ for some $x\in X$)] If $x$ is isolated, then the argument of  the subcase of isolated $x$ of the previous case applies. We assume now that $x$ is not isolated. Then, $q=1_y$ for some $y\in X$. The set $U = \{h\in Hom(XI): h(p)\in (0.5_y, 1_y]\}$ is a s desired.
	\item[\it Case \rm ($p=r_x$ for some $x\in X$ and $r\in (0,1)$)] Then $q=t_y$ for some $t\in (0,1)$ and $y\in X$. Select $a_x, a'_x, b'_x, b_x\in (0_x, 1_x)$ such that $(a_x, b_x)\subset O_p$ and $a_x<a'_x<r_x<b'_x<b_x$. Since $f$ is a homeomorphism, $f(a_x)$ and $f(b_x)$ are elements of $(0_y,1_y)$ on the opposite sides of $r_y$. The set 
$U=\{h\in Hom(XI): h(a_x')\in f((a_x, r_x)), f(b_x')\in f((r_x, b_x))\}$ is an open neighborhood of $f$. To show that $U$ is as desired, fix $h\in U$. Then $h(a_x')$ and $h(b_x')$ are on the opposite sides of $r_y$. Therefore, $h^{-1}(r_y)\in (a_x', b_x')\subset (a_x,b_x)\subset O_p$.
\end{description}
\end{proof}

\par\bigskip\noindent
\begin{lem}\label{lem:yes00}
Let $X$ be a topological space and let $XI$  be its Alexandroff Noodles. Then,  $\langle f, g\rangle \mapsto f\circ g$ is a continuous map 
from  $Hom_p(XI)\times Hom_p(XI)$ to $Hom_p(XI)$.
\end{lem}
\begin{proof}
Fix $f,g\in Hom(XI)$, $a\in XI$, and an open neighborhood  $O_c$ of $c=f(g(a))$. Put $U_{f\circ g} = \{h\in Hom(XI): h(a)\in O_c\}$. We need to find open neighborhoods $U_f$ and $U_g$ of $f$ and $g$, respectively, such  that $U_f\circ U_g \subset U_{f\circ g}$. Put $b=g(a)$. Using our agreed notations, $a\in [0_x,1_x], b\in [0_y,1_y]$, and $ c\in [0_z,1_z]$ for some $x,y,z\in X$.
\begin{description}
	\item[\it Case \rm ($a=0_x$ for some $x\in X$)] If $x$ is not isolated, then $b= 0_y$ and $c= 0_z$ for some $y,z\in X$. Put 
$U_g=\{h\in Hom(XI): h(1_x) \in (0_y, 1_y]\}$ and $U_f=\{h\in Hom(XI): h(1_y) \in (0_z, 1_z]\}$. If $g'\in U_g$ and $f'\in U_f$, then $g(0_x)=0_y$ and $f(0_y) = 0_z$. Hence, $ f'\circ g'\in U_f\circ U_g$.

\par\smallskip\noindent
If $x$ is isolated, then $b\in \{0_y,1_y\}$ and $c\in \{0_z, 1_z\}$.  All variations are treated similarly. We assume that $b=1_y$ and $c=0_z$. Put 
$U_g=\{h\in Hom(XI): h(1_x)\in [0_y,1_y)\}$ and $U_f=\{h\in Hom(XI): h(0_y)\in (0_z,1_z]\}$. To show that $U_f$ and $U_g$ are as desired, pick $g'\in U_g$ and $f'\in U_f$. Since $g'(1_x)\in [0_y,1_y)$ we conclude that $g'(1_x) = 0_y$. Hence, $g'(0_x)=1_y$. Similarly, we conclude that $f'(1_y)=0_z$. Therefore, $f'(g'(0_x))=0_z=c\in O_c$
	\item[\it Case \rm ($a=1_x$ for some $x\in X$)] If $x$ is not isolated,  then neither is  $y$ nor is  $z$. Therefore, $b= 1_y$ and $c=1_z$. Put
$U_g=\{h\in Hom(XI): h(1_x) \in (0.5_y, 1_y]\}$ and $U_f=\{h\in Hom(XI): h(1_y) \in (0.5_z, 1_z]\}$.

\par\smallskip\noindent
If $x$ is isolated, then the case is treated similarly to the case when  $a=0_x$ and $x$ is isolated.
	\item[\it Case \rm ($a=r_x$ for some $x\in X$ and $r\in (0,1)$)] Then $b= s_y$ and $c= t_z$ for some $s,t\in (0,1)$. We may assume now that $O_c=(t'_z, t''_z)\subset (0_z,1_z)$.  Pick $(s',s'')\subset (0,1)$ that contains $s$ such that $f((s_y',s_y''))\subset O_c$. Since $f$ is an automorphism, $f(s'_y)$ and $f(s''_y)$ are on the opposite sides of $t_z$. We may assume that $f(s'_y)<t_z$.
Put $U_f=\{h\in Hom(XI): h(s'_y)\in (t'_z,t_z), h(s_y'')\in(t_z,t''_z)\}$. Clearly, $f\in U_f$.  Note that any $h\in U_f$ maps $(s'_y,s''_y)$ inside of $O_c$.
\par\smallskip\noindent
Next,  pick $(r',r'')\subset (0,1)$ that contains $r$ such that $f((r_x',r_x''))\subset (s_y',s_y'')$. Since $f$ is an automorphism, $f(r'_x)$ and $f(r''_x)$ are on the opposite sides of $s_y$. We may assume that $f(r'_x)<s_y$.
Put $U_g=\{h\in Hom(XI): h(r'_x)\in (s'_y,s_y), h(r_x'')\in(s_y,s''_y)\}$. Note that any $h\in U_f$ maps $(r'_x,r''_x)$ inside of  $(s'_y,s''_y)$, which implies that $U_g$ and $U_f$ are as desired.
\end{description}
The proof is complete.
\end{proof}

\par\bigskip\noindent
Lemmas \ref{lem:yes0} and \ref{lem:yes00} imply the following.
\par\bigskip\noindent
\begin{cor}\label{cor:yes000}
Let $X$ be a topological space and let $XI$ be its Alexandroff Noodles. Then,  $Hom_p(XI)$ is a topological group.
\end{cor}

\par\bigskip\noindent
For the purpose of our next discussion, by $Hom_d(X)$ we denote the space with the underlying set $Hom(X)$ and the topology generated by sets in form $U(x,y)= \{h\in Hom(X): h(x) = y\}$. Note that the topology of $Hom_p(X)$ is a subset of the topology of $Hom_d(X)$ since 
$\{h\in Hom(X): h(x) \in O\}$ is equal to $\bigcup_{y\in O}\{h\in Hom(X): h(x) = y\}$. While $Hom_p(X)$ need not be a topological group, 
$Hom_d(X)$ is always one. To see why $f\mapsto f^{-1}$ is continuous, fix $f\in Hom(X)$ and $x\in X$.  Put $y=f(x)$ and 
$V_{f^{-1}}=\{h\in Hom(X): h(y) = x\}$. Clearly, $V_f = \{h\in Hom(X): h(x)=y\}$ is an open neighborhood of $f$ and 
$(V_f)^{-1}=V_{f^{-1}}$. A similar argument shows that the function  composition is also a continuous operation. We will use 
$Hom_d(X)$ to present our next positive observation of this study. 

\par\bigskip\noindent
Let $X$ be a space and $x, y\in X$. We say that $x$ is equivalent to $y$ if $f(x)=y$ for some $f\in Hom(X)$.
We say that $X$ is locally unique at $x$ if there exists  an open neighborhood $O$ of $x$ such that $x$ is not equivalent  to any $y$ in $O\setminus \{x\}$.  Loosely speaking, all look-alikes of $x$ are very far from $x$. There  are many such spaces among classical examples. In particular, any subspace of an ordinal is such. Moreover, any subspace of $\alpha^n$ is such for any ordinal $\alpha$ and any positive integer $n$. It is also observed in \cite{G} that scattered spaces are such as well. Observe that if $x\in X$ is not equivalent to any element in $O\setminus \{x\}$, then  $f(x)$ is not equivalent to any element in $f(O)\setminus \{f(x)\}$ for any $f\in Hom(X)$.

\par\bigskip\noindent
\begin{lem}\label{lem:yes1}
Let $X$ be locally unique at all its points. Then $Hom_p(X)$ is a topological group. Moreover,  $Hom_p(X)=Hom_d(X)$.
\end{lem}
\begin{proof}Let $\mathcal T_p$ and $\mathcal T_d$ be the topologies of $Hom_p(X)$ and $Hom_d(X)$, respectively. By our earlier discussion, it suffices to show that $\mathcal T_d\subset \mathcal T_p$. For this, fix any open neighborhood $U$ of $f$ in $Hom_p(X)$. There exist $x_1,...,x_n\in X$ and  open sets $O_1,...,O_n\subset X$ such that $f\in V= \{h\in Hom(X): h(x_i)\in O_i, i=1,...,n\}\subset U$. We may assume that each $O_i$ is small enough so that $f(x_i)$ is not equivalent to any element in $O_i\setminus \{f(x_i)\}$. Therefore, $h(x_i) = f(x_i)$ for each $i=1,...,n$ and each $h\in V$. Hence $V = \{h\in Hom(X): h(x_i) = f(x_i), i=1,...,n\}$, which completes the proof.
\end{proof}

\par\bigskip\noindent
In \cite{G}, the author proves that $Hom_p(X)=Hom_d(X)$ for scattered spaces and concludes  that $Hom_p(X)$ is a topological group. Even though Lemma \ref{lem:yes1} does not add any new interesting spaces to our collection, we included it to double down on our original claim that we have to sacrifice strong homogeneity properties in search for spaces $X$ for which $Hom_p(X)$ is a topological group.

\par\bigskip\noindent
For our next discussion, by $XX'$ we denote the Alexandroff double of $X$ . If $Y\subset X$, then $XY'$ is the corresponding subspace of the Alexandroff double of $X$. Recall that  the topology of $XX'$ is generated by sets $UU'\setminus F'$ and $F'$, where $U$ is open in $X$ and $F'$ is a finite subset of $X'$ (see \cite{Eng} for general properties of the structure).

\par\bigskip\noindent
\begin{thm}\label{thm:yes3} Let $X$ have no isolated points and
let $Hom_p(X)$ be a topological group. Then $Hom_p(XY')$ is a  topological group for any $Y\subset X$.
\end{thm}
\begin{proof} First, since $X$ has no isolated points we conclude that  $f(X)=X$ for any $f\in Hom(XY')$. In other words, $f|_X\in Hom(X)$ for every $f\in Hom(XY')$.
To show continuity of the operation of taking the inverse, fix $f\in Hom(XY')$, $p\in XY'$, and an open neighborhood $O$ of $p$. Put $q = f(p)$ and  $V = \{h\in Hom(XY'):h(q) \in O\}$. The set is an open neighborhood of $f^{-1}$ The goal is to show that there exists an open neighborhood $U$ of $f$ such that $U^{-1}\subset V$. We have two cases.

\begin{description}
	\item[\rm Case ($p$ is isolated)] Then $\{q\}$ is an open set and  $U = \{h\in Hom(XY'): h(p) = \{q\}\}$ is as desired.
	\item[\rm Case ($p$ is not isolated)] Then $p\in X$. We may assume that $O = (O_p \cup (O_p' \cap Y'))\setminus \{p'\}$ for some open neighborhood $O_p$ of $p$ in $X$. Put $V_X = \{h\in Hom(X): h(y)\in O_p\}$. Clearly, $V_X$ is an open neighborhood of $(f|_X)^{-1}$. Since $Hom_p(X)$
 is a topological group, there exists $U_X$ an open neighborhood of $f|_X$ such that $U_X^{-1}\subset V_X$. We may assume that there exist 
$x_1, x_2, ..., x_n$ and open neighborhoods $B_1,...,B_n$ of $x_1,...,x_n$ in $X$ such that 
$U_X = \{h\in Hom(X): h(x)\in O_p, h(x_i)\in B_i, \  i=1,...,n\}$.  For each $i=1,...,n$ let $O_i = B_i\cup (B_i'\cap Y')$. Put 
$U = \{h\in Hom(XY'): h(p)\in O, h(x_i)\in O_i, i=1,...,n\}$. Clearly, $U$ is an open neighborhood of $f$ in $Hom_p(XY')$ and $U^{-1}\subset V$.
\end{description}

\end{proof}

\par\bigskip\noindent
It would be interesting to know, of course, if we can drop in Theorem \ref{thm:yes3} the requirement of $X$ having non isolated points.

\par\bigskip\noindent
Theorem \ref{thm:yes3} gives us many examples. Let us isolate one very notable into a corollary.

\par\bigskip\noindent
\begin{cor}\label{cor:yes4}
$Hom_p(\mathbb R\mathbb R')$  is a topological group.
\end{cor}

\par\bigskip\noindent
We would like to finish our study with a few questions of similar character that may lead to a discovery of larger classes of spaces for which the space of automorphisms
in the topology of point-wise convergence is a topological group.

\par\bigskip\noindent
\begin{que}
Let $X$ have a basis $\mathcal B$ such that $Hom_p(B)$ is a topological group for every $B\in \mathcal B$. Is $Hom_p(X)$ a topological group?
\end{que}

\par\bigskip\noindent
\begin{que}
Let $X$ have an open cover $\mathcal U$ such that $Hom_p(U)$ is a topological group for every $U\in \mathcal U$. Is $Hom_p(X)$ a topological group?
\end{que}

\par\bigskip\noindent
\begin{que}
Let $X$ have a locally finite closed  cover $\mathcal C$ such that $Hom_p(C)$ is a topological group for every $C\in \mathcal C$. Is $Hom_p(X)$ a topological group?
\end{que}

\end{document}